\documentclass[12pt,leqno]{article}
\usepackage{amsmath, amssymb}

\newtheorem{theorem}{\bf Theorem}
\newtheorem{proposition}[theorem]{\bf Proposition}
\newtheorem{lemma}[theorem]{\bf Lemma}

\newtheorem{hypotheses}[theorem]{\bf Hypotheses}

\newcommand{\sect}[1]{\setcounter{equation}{0}\section{#1}}

\def\epsilon{\varepsilon}

\begin{document}

\Large \noindent 
{\bf Existence and nonexistence of traveling waves 
for a nonlocal monostable equation}


\vspace*{0.8em}

\normalsize
\noindent Hiroki Yagisita

\noindent
Department of Mathematics, 
Faculty of Science, 
Kyoto Sangyo University

\noindent
Motoyama, Kamigamo, Kita-Ku, Kyoto-City, 603-8555, Japan

\vfill

\noindent {\bf Abstract} \ 
We consider the nonlocal analogue of the Fisher-KPP equation 
\[u_t=\mu*u-u+f(u),\]
where $\mu$ is a Borel-measure on $\mathbb R$ with $\mu(\mathbb R)=1$ 
and $f$ satisfies $f(0)=f(1)=0$ and $f>0$ in $(0,1)$. 
We do not assume that $\mu$ is absolutely continuous 
with respect to the Lebesgue measure. 
We show that there is a constant $c_*$ 
such that it has a traveling wave solution with speed $c$ 
when $c\geq c_*$ 
while no traveling wave solution with speed $c$ when $c<c_*$, 
provided $\int_{y\in\mathbb R}e^{-\lambda y}d\mu(y)<+\infty$ 
for some positive constant $\lambda$. 
We also show that it has no traveling wave solution, 
provided 
$f^\prime(0)>0$
and 
$\int_{y\in\mathbb R}e^{-\lambda y}d\mu(y)=+\infty$ 
for all positive constants $\lambda$. 

\vspace*{0.4em} 

\noindent Keywords: spreading speed, convolution model, 

\noindent
integro-differential equation, discrete monostable equation, 

\noindent
nonlocal evolution equation, Fisher-Kolmogorov equation. 

\vspace*{2.0em}

\noindent
AMS Subject Classification: 35K57, 35K65, 35K90, 45J05.

\vspace*{1.6em}

\noindent 
A proposed running title: Traveling Waves in Nonlocal Systems II.  
                           
\newpage




\sect{Introduction}
\( \, \, \, \, \, \, \, \) 
In 1930, Fisher [8] introduced the reaction-diffusion equation $u_t=u_{xx}+u(1-u)$ 
as a model for the spatial spread of an advantageous form of a single gene 
in a population. He [9] found that there is a constant $c_*$ 
such that the equation has a traveling wave solution with speed $c$ when $c\geq c_*$ 
while it has no such solution when $c<c_*$. Kolmogorov, Petrovsky and Piskunov [16] 
obtained the same conclusion for a monostable equation $u_t=u_{xx}+f(u)$ 
with a more general nonlinearity $f$, 
and investigated long-time behavior in the model. 
Since the pioneering works, there have been extensive studies 
on traveling waves and long-time behavior for monostable evolution systems.

In this paper, we consider the following nonlocal analogue of the Fisher-KPP equation: 
\begin{equation}u_t=\mu*u-u+f(u).\end{equation} 
Here, $\mu$ is a Borel-measure on $\mathbb R$ with $\mu(\mathbb R)=1$ and the convolution 
is defined by 
\[(\mu*u)(x):=\int_{y\in\mathbb R}u(x-y)d\mu(y)\] 
for a bounded and Borel-measurable function $u$ on $\mathbb R$. The nonlinearity $f$ 
is a Lipschitz continuous function on $\mathbb R$ with $f(0)=f(1)=0$ and $f>0$ in $(0,1)$. 
Then, $G(u):=\mu*u-u+f(u)$ is a map from the Banach space $L^\infty(\mathbb R)$ 
into $L^\infty(\mathbb R)$ and it is Lipschitz continuous. 
(We note that $u(x-y)$ is a Borel-measurable function on $\mathbb R^2$, 
and $\|u\|_{L^\infty(\mathbb R)}=0$ implies 
$\|\mu*u\|_{L^1(\mathbb R)}\leq 
\int_{y\in\mathbb R}(\int_{x\in\mathbb R}|u(x-y)|dx)d\mu(y)$=0.) 
So, because the standard theory of ordinary differential equations works, 
we have well-posedness of the equation (1.1) 
and it generates a flow in $L^\infty(\mathbb R)$. 

For the nonlocal monostable equation, Atkinson and Reuter [1] 
first studied existence and nonexistence of traveling wave solutions. 
Schumacher [21, 22] showed that there is the minimal speed $c_*$ 
of traveling wave solutions 
and it has a traveling wave solution with speed $c$ when $c\geq c_*$, 
provided the extra condition 
$f(u)\leq f^\prime(0)u$ 
and some little ones. 
Here, we say that the solution $u(t,x)$ is 
{\it a traveling wave solution with profile $\psi$ and speed $c$}, 
if $u(t,x)\equiv\psi(x-x_0+ct)$ holds for some constant $x_0$ 
with $0\leq \psi\leq 1$, $\psi(-\infty)=0$ and $\psi(+\infty)=1$. 
Further, Coville, D\'{a}vila and Mart\'{i}nez [6] proved the following theorem:

\newpage

\noindent
{\bf Theorem} ([6]) \ 
{\it Suppose the nonlinearity $f\in C^1(\mathbb R)$ satisfies $f^\prime(1)<0$ 
and the Borel-measure $\mu$ has a density function $J\in C(\mathbb R)$ with 
\[\int_{y\in \mathbb R}(|y|+e^{-\lambda y})J(y)dy<+\infty\] 
for some positive constant $\lambda$. Then, there exists a constant $c_*$ 
such that the equation (1.1) has a traveling wave solution 
with monotone profile and speed $c$ when $c\geq c_*$ 
while it has no such solution when $c<c_*$.}

\vspace*{0.4em}

\noindent 
Recently, the author [28] also obtained the following: 

\vspace*{0.4em} 

\noindent
{\bf Theorem} ([28]) \ 
{\it Suppose there exists a positive constant $\lambda$ such that  
\[\int_{y\in \mathbb R}e^{\lambda|y|}d\mu(y)<+\infty\] 
holds. Then, there exists a constant $c_*$ 
such that the equation (1.1) has a traveling wave solution 
with monotone profile and speed $c$ when $c\geq c_*$ 
while it has no periodic traveling wave solution with average speed $c$ when $c<c_*$. 
Here, a solution $\{u(t,x)\}_{t\in\mathbb R}\subset L^\infty(\mathbb R)$ 
to (1.1) is said to be a periodic traveling wave solution with average speed $c$, 
if there exists a positive constant $\tau$ such that 
$u(t+\tau,x)=u(t,x+c\tau)$ holds for all $t$ and $x\in\mathbb R$ 
with $0\leq u(t,x)\leq 1$, $\lim_{x\rightarrow+\infty}u(t,x)=1$ 
and $\|u(t,x)-1\|_{L^\infty(\mathbb R)}\not=0$.} 

\vspace*{0.4em}

The goal of this paper is to improve this result of [28], 
and the following two theorems are the main results: 
\begin{theorem} \ 
Suppose there exists a positive constant $\lambda$ such that  
\[\int_{y\in \mathbb R}e^{-\lambda y}d\mu(y)<+\infty\] 
holds. Then, there exists a constant $c_*$ 
such that the equation (1.1) has a traveling wave solution 
with monotone profile and speed $c$ when $c\geq c_*$ 
while it has no periodic traveling wave solution with average speed $c$ when $c<c_*$. 
Here, a solution $\{u(t,x)\}_{t\in\mathbb R}\subset L^\infty(\mathbb R)$ 
to (1.1) is said to be a periodic traveling wave solution with average speed $c$, 
if there exists a positive constant $\tau$ such that 
$u(t+\tau,x)=u(t,x+c\tau)$ holds for all $t$ and $x\in\mathbb R$ 
with $0\leq u(t,x)\leq 1$, $\lim_{x\rightarrow+\infty}u(t,x)=1$ 
and $\|u(t,x)-1\|_{L^\infty(\mathbb R)}\not=0$.
\end{theorem}
\begin{theorem} \ 
Suppose the nonlinearity $f\in C^1(\mathbb R)$ satisfies 
\[f^\prime(0)>0.\]
Suppose the measure $\mu$ satisfies  
\[\int_{y\in \mathbb R}e^{-\lambda y}d\mu(y)=+\infty\] 
for all positive constants $\lambda$. 
Then, the equation (1.1) has no periodic traveling wave solution. 
\end{theorem}
\noindent 
In these results, we do not assume that the measure $\mu$ is absolutely continuous 
with respect to the Lebesgue measure. For example, not only the integro-differential equation 
\[\frac{\partial u}{\partial t}(t,x)=\int_0^1u(t,x-y)dy-u(t,x)+f(u(t,x))\]
but also the discrete equation 
\[\frac{\partial u}{\partial t}(t,x)=u(t,x-1)-u(t,x)+f(u(t,x))\]
satisfies the assumption of Theorem 1 for the measure $\mu$. 
In order to prove these results, 
we employ the recursive method for monotone dynamical systems 
by Weinberger [25] and Li, Weinberger and Lewis [17]. 
We note that the semiflow generated by the equation (1.1) 
does not have compactness with respect to the compact-open topology. 

Schumacher [21, 22], Carr and Chmaj [3] and Coville, D\'{a}vila and Mart\'{i}nez [6] 
also studied uniqueness of traveling wave solutions. 
In [6], we could see an interesting example of nonuniqueness, 
where the equation (1.1) admits infinitely many monotone profiles 
for standing wave solutions but it admits no continuous one. 
See, e.g., [5, 7, 10, 11, 12, 13, 14, 15, 18, 19, 23, 24, 26, 27] 
on traveling waves and long-time behavior in various monostable evolution systems, 
[2, 4] nonlocal bistable equations and [20] Euler equation. 

In Section 2, we recall abstract results for monotone semiflows from [28]. 
In Section 3, we give basic facts for nonlocal equations in $L^\infty(\mathbb R)$. 
In Section 4, we prove Theorem 1. 
In Section 5, we recall a result on spreading speeds by Weinberger [25].  
In Section 6, we prove Theorem 2.  
  
\sect{Abstract results for monotone semiflows}
\( \, \, \, \, \, \, \, \) 
In this section, we recall some abstract results for monotone semiflows 
from [28]. Put a set of functions on $\mathbb R$; 
\[\mathcal M:=\{u\, |\, u \text{ is a monotone nondecreasing}\]
\[\text{ and left continuous function on } \mathbb R \text{ with } 0\leq u\leq 1\}.\] 
The followings are basic conditions for discrete dynamical systems on $\mathcal M$: 
\begin{hypotheses} \ 
Let $Q_0$ be a map from $\mathcal M$ into $\mathcal M$. 

{\rm (i)} \ $Q_0$ is continuous in the following sense: 
If a sequence $\{u_k\}_{k\in \mathbb N}\subset \mathcal M$ 
converges to $u\in \mathcal M$ uniformly on every bounded interval, 
then the sequence $\{Q_0[u_k]\}_{k\in \mathbb N}$ converges to $Q_0[u]$ almost everywhere. 

{\rm (ii)} \ $Q_0$ is order preserving; i.e., 
\[u_1\leq u_2 \Longrightarrow Q_0[u_1]\leq Q_0[u_2]\] 
for all $u_1$ and $u_2\in \mathcal M$. 
Here, $u\leq v$ means that $u(x)\leq v(x)$ holds for all $x\in \mathbb R$. 

{\rm (iii)} \ $Q_0$ is translation invariant; i.e., 
\[T_{x_0}Q_0=Q_0T_{x_0}\] 
for all $x_0\in \mathbb R$. Here, $T_{x_0}$ is the translation operator 
defined by $(T_{x_0}[u])(\cdot):=u(\cdot-x_0)$. 

{\rm (iv)} \ $Q_0$ is monostable; i.e., 
\[0<\alpha<1 \Longrightarrow \alpha < Q_0[\alpha]\] 
for all constant functions $\alpha$. 
\end{hypotheses} 
{\bf Remark} \ If $Q_0$ satisfies Hypothesis 3 (iii), 
then $Q_0$ maps constant functions to constant functions. 

\vspace*{0.4em} 

\noindent 
We add the following conditions to Hypotheses 3 
for continuous dynamical systems on $\mathcal M$: 
\begin{hypotheses} \ 
Let $Q:=\{Q^t\}_{t\in[0,+\infty)}$ be a family of maps 
from $\mathcal M$ to $\mathcal M$. 

{\rm (i)} \ $Q$ is a semigroup; 
i.e., $Q^t\circ Q^s=Q^{t+s}$ for all $t$ and $s\in [0,+\infty)$. 

{\rm (ii)} \ $Q$ is continuous in the following sense: 
Suppose a sequence $\{t_k\}_{k\in \mathbb N}\subset[0,+\infty)$ converges to $0$, 
and $u\in \mathcal M$. Then, the sequence $\{Q^{t_k}[u]\}_{k\in \mathbb N}$ 
converges to $u$ almost everywhere. 
\end{hypotheses}

From [28], we recall the following two results
for continuous dynamical systems on $\mathcal M$: 
\begin{theorem} \ 
Let $Q^t$ be a map from $\mathcal M$ to $\mathcal M$ for $t\in[0,+\infty)$. 
Suppose $Q^t$ satisfies {\rm Hypotheses 3} for all $t\in (0,+\infty)$, 
and $Q:=\{Q^t\}_{t\in[0,+\infty)}$ {\rm Hypotheses 4}. Then, the following holds {\rm :} 

Let $c\in\mathbb R$. Suppose there exist $\tau\in(0,+\infty)$ and $\phi\in\mathcal M$ 
with $(Q^\tau[\phi])(x-c\tau)\leq \phi(x)$, $\phi\not\equiv 0$ and $\phi\not\equiv 1$. 
Then, there exists $\psi\in\mathcal M$ with $\psi(-\infty)=0$ and $\psi(+\infty)=1$ 
such that $(Q^t[\psi])(x-ct)\equiv \psi(x)$ holds for all $t\in[0,+\infty)$. 
\end{theorem}
\begin{theorem} \ 
Let $Q^t$ be a map from $\mathcal M$ to $\mathcal M$ for $t\in[0,+\infty)$. 
Suppose $Q^t$ satisfies {\rm Hypotheses 3} for all $t\in (0,+\infty)$, 
and $Q:=\{Q^t\}_{t\in[0,+\infty)}$ {\rm Hypotheses 4}. 
Then, there exists $c_*\in (-\infty, +\infty]$ 
such that the following holds {\rm :} 

Let $c\in \mathbb R$. Then, there exists $\psi\in\mathcal M$ 
with $\psi(-\infty)=0$ and $\psi(+\infty)=1$ 
such that $(Q^t[\psi])(x-ct)\equiv \psi(x)$ holds for all $t\in[0,+\infty)$ 
if and only if $c\geq c_*$. 
\end{theorem}

\sect{Basic facts for nonlocal equations in $L^\infty(\mathbb R)$} 
\( \, \, \, \, \, \, \, \) 
In this section, we give some basic facts for the equation 
\begin{equation}u_t=\hat\mu*u+g(u)\end{equation} 
on the phase space $L^\infty(\mathbb R)$. 
First, we have the comparison theorem for (3.1) on $L^\infty(\mathbb R)$: 
\begin{lemma} \ 
Let $\hat \mu$ be a Borel-measure on $\mathbb R$ with $\hat \mu (\mathbb R)<+\infty$.
Let $g$ be a Lipschitz continuous function on $\mathbb R$. 
Let $T\in(0,+\infty)$, and two functions $u^1$ 
and $u^2\in C^1([0,T],L^\infty(\mathbb R))$. 
Suppose that for any $t\in[0,T]$, 
the inequality 
\[u^1_t-\left(\hat\mu*u^1+g(u^1)\right) 
\leq u^2_t-\left(\hat\mu*u^2+g(u^2)\right)\] 
holds almost everywhere in $x$. Then, 
the inequality $u^1(T,x)\leq u^2(T,x)$ holds almost everywhere in $x$ 
if the inequality $u^1(0,x)\leq u^2(0,x)$ holds almost everywhere in $x$. 
\end{lemma}
{\bf Proof.} \ 
Put $K\in\mathbb R$ by 
\begin{equation}K:=-\inf_{h>0,u\in\mathbb R}\frac{g(u+h)-g(u)}{h},\end{equation} 
and $v\in C^1([0,T],L^\infty(\mathbb R))$ by 
\begin{equation}v(t):=e^{Kt}(u^2-u^1)(t).\end{equation} 
Then, we have the ordinary differential equation 
\begin{equation}\frac{dv}{dt}=F(t,v)\end{equation} 
in $L^\infty(\mathbb R)$ with $v(0)=(u^2-u^1)(0)$ 
as we define a map 
$F:[0,T]\times L^\infty(\mathbb R)\rightarrow L^\infty(\mathbb R)$ by 
\[F(t,w):=\hat \mu*w+Kw
+e^{Kt}\left(g(u^1(t)+e^{-Kt}w)-g(u^1(t))\right)
+e^{Kt}a(t),\] 
where 
\[a:=\left(\frac{du^2}{dt}-\left(\hat \mu*u^2+g(u^2)\right)\right)
-\left(\frac{du^1}{dt}-\left(\hat \mu*u^1+g(u^1)\right)\right).\] 
For any $t\in[0,T]$, we see the inequality 
\begin{equation}a(t,x)\geq 0\end{equation}
almost everywhere in $x$.
Take the solution $\tilde{v}\in C^1([0,T],L^\infty(\mathbb R))$ to 
\begin{equation}\tilde{v}(t)=v(0)+\int_0^t\max\{F(s,\tilde{v}(s)),0\}ds.\end{equation} 
Then, for any $t\in[0,T]$, we have 
\begin{equation}\tilde{v}(t,x)\geq v(0,x)=(u^2-u^1)(0,x)\geq 0\end{equation}
almost everywhere in $x$. By using (3.2), (3.5) and (3.7), for any $t\in[0,T]$, 
we also have the inequality $F(t,\tilde{v}(t))\geq 0$ almost everywhere in $x$. 
Hence, from (3.6), $\tilde{v}(t)$ is the solution 
to the same ordinary differential equation (3.4) in $L^\infty(\mathbb R)$ 
as $v(t)$ with $\tilde{v}(0)=v(0)$. So, in virtue of (3.3) and (3.7), 
\[(u^2-u^1)(T,x)=e^{-KT}v(T,x)=e^{-KT}\tilde{v}(T,x)\geq 0\]
holds almost everywhere in $x$. 
\hfill 
$\blacksquare$ 

\vspace*{0.8em} 

The following lemma gives a invariant set and some positively invariant sets 
of the flow on $L^\infty(\mathbb R)$ 
generated by the equation (3.1): 
\begin{lemma} \ 
Let $\hat\mu$ be a Borel-measure on $\mathbb R$ with $\hat\mu(\mathbb R)<+\infty$. 
Let $g$ be a Lipschitz continuous function on $\mathbb R$. 
Then, the followings hold {\rm :}

{\rm (i)} \ 
For any $u_0\in BC(\mathbb R)$, 
there exists a solution $\{u(t)\}_{t\in\mathbb R}\subset BC(\mathbb R)$ 
to (3.1) with $u(0)=u_0$. Here, $BC(\mathbb R)$ denote the set 
of bounded and continuous functions on $\mathbb R$. 

{\rm (ii)} \ Suppose a constant $\gamma$ 
satisfies $\gamma\hat\mu(\mathbb R)+g(\gamma)=0$.  
If $u_0\in L^\infty(\mathbb R)$ satisfies $\gamma \leq u_0$, then 
there exists a solution $\{u(t)\}_{t\in[0,+\infty)}\subset L^\infty(\mathbb R)$ 
to (3.1) with $u(0)=u_0$ and $\gamma\leq u(t)$. 
If $u_0\in L^\infty(\mathbb R)$ satisfies $u_0\leq \gamma$, then 
there exists a solution $\{u(t)\}_{t\in[0,+\infty)}\subset L^\infty(\mathbb R)$ 
to (3.1) with $u(0)=u_0$ and $u(t)\leq \gamma$. 

{\rm (iii)} \ 
If $u_0$ is a bounded and monotone nondecreasing function on $\mathbb R$, 
then there exists a solution $\{u(t)\}_{t\in[0,+\infty)}\subset L^\infty(\mathbb R)$ 
to (3.1) with $u(0)=u_0$ such that $u(t)$ is a bounded 
and monotone nondecreasing function on $\mathbb R$ for all $t\in[0,+\infty)$. 
If $u_0$ is a bounded and monotone nonincreasing function on $\mathbb R$, 
then there exists a solution $\{u(t)\}_{t\in[0,+\infty)}\subset L^\infty(\mathbb R)$ 
to (3.1) with $u(0)=u_0$ such that $u(t)$ is a bounded 
and monotone nonincreasing function on $\mathbb R$ for all $t\in[0,+\infty)$. 
\end{lemma}
{\bf Proof.} \ We could see (i), because 
$BC(\mathbb R)$ is a closed sub-space of the Banach space $L^\infty(\mathbb R)$ 
and $u\in BC(\mathbb R)$ implies $\hat\mu*u+g(u)\in BC(\mathbb R)$. 

We could also see (ii) by using Lemma 7, because 
the constant $\gamma$ is a solution to (3.1).  
 
We show (iii). Suppose $u_0$ is a bounded and monotone nondecreasing 
function on $\mathbb R$. 
We take a solution $\{u(t)\}_{t\in[0,+\infty)}\subset L^\infty(\mathbb R)$ 
to (3.1) with $u(0)=u_0$. Let $t\in[0,+\infty)$ and $h\in[0,+\infty)$. 
Then, by Lemma 7, we see $u(t,x)\leq u(t,x+h)$ almost everywhere in $x$. 
We take a cutoff function $\rho\in C^\infty(\mathbb R)$ with 
\[|x|\geq 1/2 \, \Longrightarrow \, \rho(x)=0,\]
\[|x|<1/2 \, \Longrightarrow \, \rho(x)>0\] 
and
\[\int_{x\in\mathbb R}\rho(x)dx=1.\] 
As we put  
\[v_n(x):=\int_{y\in \mathbb R}2^n\rho(2^n(x-y))u(t,y)dy\]
for $n\in\mathbb N$, we see $v_n(x)\leq v_n(x+h)$ for all $x\in\mathbb R$. 
Therefore, $v_n$ is smooth, bounded and monotone nondecreasing. 
By Helly's theorem, there exist a subsequence $n_k$ 
and a bounded and monotone nondecreasing function $\psi$ on $\mathbb R$ 
such that $\lim_{k\rightarrow\infty}v_{n_k}(x)=\psi(x)$ holds for all $x\in \mathbb R$. 
Then, $\|u(t,x)-\psi(x)\|_{L^1([-C,+C])}
\leq \lim_{k\rightarrow\infty}
(\|u(t,x)-v_{n_k}(x)\|_{L^1([-C,+C])}
+\|v_{n_k}(x)-\psi(x)\|_{L^1([-C,+C])})=0$ holds for all $C\in (0,+\infty)$. 
Hence, we obtain $\|u(t,x)-\psi(x)\|_{L^\infty(\mathbb R)}=0$. 
\hfill 
$\blacksquare$ 

\vspace*{0.4em} 

\begin{lemma} \ 
Let $\hat \mu$ be a Borel-measure on $\mathbb R$ with $\hat \mu (\mathbb R)<+\infty.$ 
Let $\{u_n\}_{n=1}^\infty$ be a sequence of bounded and continuous functions on $\mathbb R$ with 
\[\sup_{n\in\mathbb N, x\in\mathbb R}|u_n(x)|<+\infty.\]
Suppose the sequence $\{u_n\}_{n=1}^\infty$ converges to $0$ uniformly on every bounded interval. 
Then, the sequence $\{\hat\mu*u_n\}_{n=1}^\infty$ converges to $0$ uniformly on every bounded interval. 
\end{lemma}
{\bf Proof.} \ 
Let $\varepsilon\in(0,+\infty)$. 
We take a positive constant $C$ such that 
\[\left(\sup_{n\in\mathbb N, x\in\mathbb R}|u_n(x)|\right)
\hat\mu(\mathbb R\setminus (-C,+C))
\leq \varepsilon\]
holds. Then, because 
\[
|(\hat\mu*u_n)(x)|
\leq\int_{y\in(-C,+C)}|u_n(x-y)|d\hat\mu(y)
+\int_{y\in\mathbb R\setminus(-C,+C)}|u_n(x-y)|d\hat\mu(y)
\]
\[\leq \left(\sup_{y\in(-C,+C)}|u_n(x-y)|\right)\hat\mu(\mathbb R)
+\left(\sup_{y\in\mathbb R}|u_n(x-y)|\right)
\hat\mu(\mathbb R\setminus (-C,+C))
\]
holds, we have 
\[
\sup_{x\in[-I,+I]}|(\hat\mu*u_n)(x)|
\leq \left(\sup_{y\in(-(I+C),+(I+C))}|u_n(y)|\right)\hat\mu(\mathbb R) 
+\varepsilon
\]
for all $I\in(0,+\infty)$. Hence, we obtain 
\[\limsup_{n\rightarrow\infty}\sup_{x\in[-I,+I]}|(\hat\mu*u_n)(x)|
\leq \varepsilon\]
for all $I\in(0,+\infty)$. 
\hfill 
$\blacksquare$ 

\begin{proposition} \ 
Let $\hat \mu$ be a Borel-measure on $\mathbb R$ with $\hat \mu (\mathbb R)<+\infty$, 
$g$ a Lipschitz continuous function on $\mathbb R$, and $T$ a positive constant. 
Let a sequence $\{u_n\}_{n=0}^\infty\subset C^1([0,T],L^\infty(\mathbb R))$ 
of solutions to the equation (3.1) satisfy 
\[\sup_{n\in\mathbb N,x\in\mathbb R}|u_n(0,x)-u_0(0,x)|<+\infty.\]  
Suppose \[\lim_{n\rightarrow\infty}\sup_{x\in\mathbb [-I,+I]}|u_n(0,x)-u_0(0,x)|=0\] 
holds for all positive constants $I$. Then, 
\[\lim_{n\rightarrow\infty}\sup_{t\in[0,T]}\|u_n(t,x)-u_0(t,x)\|_{L^\infty([-J,+J])}=0\] 
holds for all positive constants $J$. 
\end{proposition} 
{\bf Proof.} \ 
First, we take a sequence $\{w_n\}_{n=1}^\infty$ 
of nonnegative, bounded and continuous functions on $\mathbb R$ with 
\begin{equation}\sup_{n\in\mathbb N, x\in\mathbb R}|w_n(x)|<+\infty\end{equation}
such that 
$\{w_n\}_{n=1}^\infty$ converges to $0$ uniformly on every bounded interval and 
\begin{equation}|u_n(0,x)-u_0(0,x)|\leq w_n(x)\end{equation}
holds for all $n\in\mathbb N$ and $x\in\mathbb R$. 
Let $\hat A$ denote the bounded and linear operator 
from the Banach space $BC(\mathbb R)$ to $BC(\mathbb R)$ 
defined by 
\[\hat Aw:=\hat\mu*w.\] 

From (3.8), we see $\sup_{n\in\mathbb N, x\in\mathbb R}|(\hat A^k w_n)(x)|<+\infty$ 
for all $k=0,1,2,\cdots$. Hence, because of 
$\lim_{n\rightarrow\infty}\sup_{x\in[-I,+I]}|w_n(x)|=0$ 
for all $I\in(0,+\infty)$, by Lemma 9, we have 
\begin{equation}\lim_{n\rightarrow\infty}
\sup_{x\in[-J,+J]}|(\hat A^k w_n)(x)|=0\end{equation} 
for all $J\in(0,+\infty)$ and $k=0,1,2,\cdots$. 

Let $\gamma$ denote the constant defined by 
\[\gamma:=\sup_{h>0,u\in\mathbb R}\frac{g(u+h)-g(u)}{h}.\] 
Then, we consider the following two sequences $\{\underline{v}_n\}_{n=1}^\infty$ 
and $\{\overline{v}_n\}_{n=1}^\infty\subset C^1
$
$
([0,T],L^\infty(\mathbb R))$ defined by  
\[\underline{v}_n(t,x):=u_0(t,x)-e^{\gamma t}(e^{\hat At}w_n)(x)\]
and
\[\overline{v}_n(t,x):=u_0(t,x)+e^{\gamma t}(e^{\hat At}w_n)(x).\]
Because  $(e^{\hat At}w_n)(x)$ is nonnegative for all $n\in\mathbb N$, 
$t\in[0,+\infty)$ and $x\in\mathbb R$, 
the function $\underline{v}_n$ is a sub-solution to (3.1) 
and $\overline{v}_n$ is a super-solution to (3.1) for all $n\in\mathbb N$. 
So, by Lemma 7 and (3.9), for any $n\in\mathbb N$ and $t\in[0,T]$, 
\begin{equation}|u_n(t,x)-u_0(t,x)|\leq e^{\gamma t}(e^{\hat At}w_n)(x)\end{equation}
holds almost everywhere in $x$. 

Let $\varepsilon\in(0,+\infty)$. We take $N\in\mathbb N$ such that 
\[(1+e^{\gamma T})
\left(\sum_{k=N}^\infty\frac{(T\|\hat A\|_{BC(\mathbb R)\rightarrow BC(\mathbb R)})^k}{k!}\right)
\left(\sup_{n\in\mathbb N, x\in\mathbb R}|w_n(x)|\right)
\leq\varepsilon\] 
holds. Then, in virtue of (3.11), we see 
\[\|u_n(t,x)-u_0(t,x)\|_{L^\infty([-J,+J])}
\leq \sup_{x\in[-J,+J]}|e^{\gamma t}(e^{\hat At}w_n)(x)|\]
\[= e^{\gamma t}\left(\sup_{x\in[-J,+J]}\left|
\left(\sum_{k=0}^{N-1}\frac{t^k}{k!}(\hat A^k w_n)(x)\right)
+\left(\left(\sum_{k=N}^\infty\frac{t^k}{k!}\hat A^k\right)w_n\right)(x)
\right|\right)\]
\[\leq (1+e^{\gamma T})\left(\sum_{k=0}^{N-1}\frac{T^k}{k!}
\left(\sup_{x\in[-J,+J]}|(\hat A^k w_n)(x)|\right)
\right)+\varepsilon\] 
for all $J\in(0,+\infty)$, $n\in\mathbb N$ and $t\in[0,T]$. So, by (3.10), we obtain 
\[\limsup_{n\rightarrow\infty}\sup_{t\in[0,T]}\|u_n(t,x)-u_0(t,x)\|_{L^\infty([-J,+J])}
\leq \varepsilon\]
for all $J\in(0,+\infty)$. 
\hfill 
$\blacksquare$ 

\sect{Proof of Theorem 1} 
\( \, \, \, \, \, \, \, \) 
In this section, we prove Theorem 1 by using the results of Sections 2 and 3.  
The argument in this section is almost similar as in [28]. 
First, we recall that 
$\mu$ is a Borel-measure on $\mathbb R$ with $\mu(\mathbb R)=1$,  
$f$ is a Lipschitz continuous function on $\mathbb R$ 
with $f(0)=f(1)=0$ and $f>0$ in $(0,1)$ 
and 
the set $\mathcal M$ has been defined at the beginning of Section 2. 
Then, in virtue of Lemmas 7, 8 and Proposition 10, $Q^t \ (t\in(0,+\infty))$ 
satisfies Hypotheses 3 and $Q$ Hypotheses 4 
for the semiflow $Q=\{Q^t\}_{t\in[0,+\infty)}$ on $\mathcal M$ 
generated by (1.1). So, Theorems 5 and 6 can work for this semiflow on $\mathcal M$. 

\vspace*{0.8em} 

If the flow on $L^\infty(\mathbb R)$ generated by (1.1) 
has a {\it periodic} traveling wave solution 
with {\it average} speed $c$ (even if the profile is not a monotone function), 
then it has a traveling wave solution with {\it monotone} profile and speed $c$:  
\begin{theorem} \ 
Let $c\in\mathbb R$. Suppose there exist a positive constant $\tau$ 
and a solution $\{u(t,x)\}_{t\in\mathbb R} 
\subset L^\infty(\mathbb R)$ to (1.1) with $0\leq u(t,x)\leq 1$, 
$\lim_{x\rightarrow+\infty} u(t,x)=1$ and $\|u(t,x)-1\|_{L^\infty(\mathbb R)}\not=0$ 
such that \[u(t+\tau,x)=u(t,x+c\tau)\] holds for all $t$ and $x\in\mathbb R$.  
Then, there exists $\psi\in\mathcal M$ with $\psi(-\infty)=0$ and $\psi(+\infty)=1$ 
such that $\{\psi(x+ct)\}_{t\in\mathbb R}$ is a solution to (1.1). 
\end{theorem}
{\bf Proof.} \ 
Put two monotone nondecreasing functions 
$\varphi(x):=\max\{\alpha\in\mathbb R \, | \, 
\alpha \leq u(0,y) 
\text{ holds almost everywhere in } y\in (x,+\infty)\}$ 
and 
$\phi(x):=\lim_{h\downarrow+0}
$
$
\varphi(x-h)$. Then, $\phi\in\mathcal M$, 
$\phi(-\infty)<1$ and $\phi(+\infty)=1$ hold.
We take a cutoff function $\rho\in C^\infty(\mathbb R)$ with 
\[|x+1/2|\geq 1/2 \, \Longrightarrow \, \rho(x)=0,\]
\[|x+1/2|<1/2 \, \Longrightarrow \, \rho(x)>0\] 
and
\[\int_{x\in\mathbb R}\rho(x)dx=1.\] 
As we put  
\[v_n(x):=\int_{y\in \mathbb R}2^n\rho(2^n(x-y))u(0,y)dy\]
for $n\in\mathbb N$, we see $\phi\leq v_n$. Let $N\in\mathbb N$. 
Because of $\lim_{n\rightarrow\infty}\|v_{n}(x)-u(0,x)\|_{L^1([-N,+N])}=0$, 
there exists a subsequence $n_k$ such that 
$\lim_{k\rightarrow\infty}v_{n_k}(x)
$ 
$
=u(0,x)$ 
almost everywhere in $x\in[-N,+N]$. Therefore, we have 
$\phi(x)\leq u(0,x)$ almost everywhere in $x\in\mathbb R$. 
So, by Lemma 7, we obtain 
$Q^\tau[\phi](x-c\tau)\leq u(\tau,x-c\tau)=u(0,x)$ almost everywhere in $x$. 
Hence, because $Q^\tau[\phi](x-c\tau)\leq \varphi(x)$ holds, 
we get $Q^\tau[\phi](x-c\tau)\leq \phi(x)$. Therefore, by Theorem 5, 
there exists $\psi\in\mathcal M$ with $\psi(-\infty)=0$ and $\psi(+\infty)=1$
such that $Q^t[\psi](x-ct)\equiv\psi(x)$ holds for all $t\in[0,+\infty)$. 
\hfill 
$\blacksquare$ 

\vspace*{0.8em} 

The infimum $c_*$ of the speeds of traveling wave solutions is not $-\infty$, 
and there is a traveling wave solution with speed $c$ when $c\geq c_*$: 
\begin{lemma} \ 
There exists $c_*\in(-\infty,+\infty]$ such that the following holds {\rm :} 

Let $c\in \mathbb R$. Then, there exists $\psi\in\mathcal M$ 
with $\psi(-\infty)=0$ and $\psi(+\infty)=1$ 
such that $\{\psi(x+ct)\}_{t\in\mathbb R}$ is a solution to (1.1) 
if and only if $c\geq c_*$. 
\end{lemma}
{\bf Proof.} \ 
It follows from Theorem 6. 
\hfill 
$\blacksquare$ 

\vspace*{0.8em} 

\noindent
{\bf Proof of Theorem 1.} 

Let $c_*$ be the infimum  of the speeds of traveling wave solutions
with monotone profile. Then, in virtue of Theorem 11 and Lemma 12, 
it is sufficient if we show $c_*\not=+\infty$. 

Take $K\in[0,+\infty)$ such that  
\[K\geq\max\left\{\int_{y\in\mathbb R}e^{-\lambda y}d\mu(y),\mu(\mathbb R)\right\}
-1+\sup_{h>0}\frac{f(h)}{h}.\] 
As we put $\phi(x):=\min\{e^{\lambda x},1\}\in\mathcal M$, 
we see 
\[(\mu*\phi)(x)\leq\min
\left\{\left(\int_{y\in\mathbb R}e^{-\lambda y}d\mu(y)\right)e^{\lambda x},
\mu(\mathbb R)\right\}\]
\[\leq \max\left\{\int_{y\in\mathbb R}e^{-\lambda y}d\mu(y),\mu(\mathbb R)\right\}\phi(x).\]
So, $e^{Kt}\phi(x)$ is a super-solution to (1.1), because of 
\[e^{Kt}(\mu*\phi)-e^{Kt}\phi+f(e^{Kt}\phi)\leq Ke^{Kt}\phi.\] 
Hence, by Lemma 7, we obtain 
$Q^1[\phi](x)\leq e^K\phi(x)\leq e^{\lambda(x+\frac{K}{\lambda})}$, 
and $Q^1[\phi](x-\frac{K}{\lambda})\leq \phi(x)$. Therefore, from Theorem 5, 
there exists $\psi\in\mathcal M$ with $\psi(-\infty)=0$ and $\psi(+\infty)=1$
such that $Q^t[\psi](x-\frac{K}{\lambda}t)\equiv\psi(x)$ holds for all $t\in[0,+\infty)$. 
So, $c_*\leq \frac{K}{\lambda}$ holds. 
\hfill 
$\blacksquare$

\sect{A result on spreading speeds by Weinberger} 
\( \, \, \, \, \, \, \, \) 
In this section, we recall a result by Weinberger [25]. 
It is used to prove Theorem 2 in Section 6. 
Put a set of functions on $\mathbb R$; 
\[\mathcal B:=\{u\, |\, u \text{ is a continuous function on } 
\mathbb R \text{ with } 0\leq u\leq 1\}.\] 
\begin{hypotheses} \ 
Let $\tilde Q_0$ be a map from $\mathcal B$ into $\mathcal B$. 

{\rm (i)} \ $\tilde Q_0$ is continuous in the following sense: 
If a sequence $\{u_k\}_{k\in \mathbb N}\subset \mathcal B$ 
converges to $u\in \mathcal B$ uniformly on every bounded interval, 
then the sequence $\{(\tilde Q_0[u_k])(x)\}_{k\in \mathbb N}$ 
converges to $(\tilde Q_0[u])(x)$ for all $x\in\mathbb R$. 

{\rm (ii)} \ $\tilde Q_0$ is order preserving; i.e., 
\[u_1\leq u_2 \Longrightarrow \tilde Q_0[u_1]\leq \tilde Q_0[u_2]\] 
for all $u_1$ and $u_2\in \mathcal B$. 
Here, $u\leq v$ means that $u(x)\leq v(x)$ holds for all $x\in \mathbb R$. 

{\rm (iii)} \ $\tilde Q_0$ is translation invariant; i.e., 
\[T_{x_0}\tilde Q_0=\tilde Q_0T_{x_0}\] 
for all $x_0\in \mathbb R$. Here, $T_{x_0}$ is the translation operator 
defined by $(T_{x_0}[u])(\cdot):=u(\cdot-x_0)$. 

{\rm (iv)} \ $\tilde Q_0$ is monostable; i.e., 
\[0<\alpha<1 \Longrightarrow \alpha < \tilde Q_0[\alpha]\] 
for all constant functions $\alpha$, and $\tilde Q_0[0]=0$.  
\end{hypotheses}
{\bf Remark} \ If $\tilde Q_0$ satisfies Hypotheses 13 (ii) and (iii), 
then $\tilde Q_0$ maps monotone functions to monotone functions. 

\begin{theorem} \ 
Let a map $\tilde Q_0 : \mathcal B \rightarrow \mathcal B$ 
satisfy Hypotheses 13. 
Let a continuous and monotone nonincreasing function $\varphi$ on $\mathbb R$ 
with $0<\varphi(-\infty)<1$ satisfy $\varphi(x)=0$ for all $x\in[0,+\infty)$. 
For $c\in \mathbb R$, define the sequence $\{a_{c,n}\}_{n=0}^\infty$ 
of continuous and monotone nonincreasing functions on $\mathbb R$ 
by the recursion 
\[a_{c,n+1}(x):=\max\{(\tilde Q_0[a_{c,n}])(x+c),\varphi(x)\}\]
with $a_{c,0}:=\varphi$. Then, 
\[0\leq a_{c,n}\leq a_{c,n+1}\leq 1\] 
holds for all $c\in\mathbb R$ and $n=0,1,2,\cdots$. 
For $c\in \mathbb R$, define the bounded 
and monotone nonincreasing function $a_{c}$ on $\mathbb R$ 
by \[a_{c}(x):=\lim_{n\rightarrow\infty}a_{c,n}(x).\] 

Let $\tilde \nu$ be a Borel-measure on $\mathbb R$ with $1<\tilde \nu(\mathbb R)<+\infty$. 
Suppose there exists a positive constant $\varepsilon$ 
such that the inequality 
\[\tilde \nu*u\leq \tilde Q_0[u]\] 
holds for all $u\in\mathcal B$ with $u\leq \varepsilon$. Then, the inequality 
\[\inf_{\lambda>0}\frac{1}{\lambda}\log\int_{y\in\mathbb R}
e^{\lambda y}d\tilde \nu(y)
\, 
\leq
\, 
\sup\{c\in\mathbb R \, | \, a_{c}(+\infty)=1\}\] 
holds. 
\end{theorem}   
{\bf Proof.} \   
It follows from Lemma 5.4 and Theorem 6.4 in [25] 
with $N:=1$, $\mathcal H:=\mathbb R$, $\pi_0:=0$, $\pi_1=\pi_+:=1$, 
$S^{N-1}:=\{\pm 1\}$ and $\xi:=+1$. 
\hfill 
$\blacksquare$ 

\vspace*{0.8em}

From Theorem 14, we have the following: 
\begin{proposition} \
Let $\hat \mu$ be a Borel-measure on $\mathbb R$ with $\hat \mu (\mathbb R)=1$. 
Let $c_0\in\mathbb R$, and $\hat \psi$ be a monotone nonincreasing function 
on $\mathbb R$ with $\hat\psi(-\infty)=1$ and $\hat\psi(+\infty)=0$. 
Suppose $\{\hat\psi (x-c_0t)\}_{t\in\mathbb R}\subset L^\infty(\mathbb R)$ 
is a solution to 
\begin{equation}u_t=\hat \mu*u-u+f(u).\end{equation} 

Let $\, \tilde Q_0:\, \mathcal B\rightarrow \mathcal B\, $ be the time $1$ map 
of the semiflow on $\mathcal B$ generated by the equation (5.1). 
Let $\tilde \nu$ be a Borel-measure on $\mathbb R$ with $1<\tilde \nu(\mathbb R)<+\infty$. 
Suppose there exists a positive constant $\varepsilon$ 
such that the inequality 
\[\tilde \nu*u\leq \tilde Q_0[u]\] 
holds for all $u\in\mathcal B$ with $u\leq \varepsilon$. Then, the inequality 
\[\inf_{\lambda>0}\frac{1}{\lambda}\log\int_{y\in\mathbb R}
e^{\lambda y}d\tilde \nu(y)
\, 
\leq
\, 
c_0\] 
holds. 
\end{proposition}
{\bf Proof.} \   
We take a continuous and monotone nonincreasing function $\varphi$ on $\mathbb R$ 
with $0<\varphi(-\infty)<1$ and $\varphi(x)=0$ for all $x\in[0,+\infty)$. 
For $c\in \mathbb R$, we define the sequence $\{a_{c,n}\}_{n=0}^\infty$ 
of continuous and monotone nonincreasing functions on $\mathbb R$ 
by the recursion 
\[a_{c,n+1}(x):=\max\{(\tilde Q_0[a_{c,n}])(x+c),\varphi(x)\}\]
with $a_{c,0}:=\varphi$. 
We also take $x_0\in\mathbb R$ such that 
\[\varphi(x)\leq \hat\psi(x-x_0)\]
holds for all $x\in\mathbb R$. 

Let $c\in[c_0,+\infty)$. Then, we show $a_{c,n}(x)\leq \hat\psi(x-x_0)$ 
for all $n=0,1,2,\cdots$. 
We have $a_{c,0}(x)=\varphi(x)\leq \hat\psi_0(x-x_0)$. 
As $a_{c,n}(x)\leq \hat\psi(x-x_0)$ holds almost everywhere in $x$, 
\[
a_{c,n+1}(x)
\leq \max\{(\tilde Q_0[a_{c,n}])(x+c_0),\varphi(x)\}
\]
\[
\leq \max\{\hat\psi(x-x_0),\varphi(x)\}
=\hat\psi(x-x_0)
\]
also holds almost everywhere in $x$, because 
$\hat\psi(x-x_0-c_0t)$ is a solution to (5.1). 
So, for any $n=0,1,2,\cdots$, the inequality 
$a_{c,n}(x)\leq \hat\psi(x-x_0)$ 
holds almost everywhere in $x$. Hence, because 
$a_{c,n}$ is continuous and $\hat\psi$ is monotone, we have
\begin{equation}a_{c,n}(x)\leq \hat\psi(x-x_0)\end{equation}
for all $x\in\mathbb R$, $c\in[c_0,+\infty)$ and $n=0,1,2,\cdots$. 
Therefore, by Theorem 14, (5.2) and $\hat\psi(+\infty)=0$, the inequality 
\[\inf_{\lambda>0}\frac{1}{\lambda}\log\int_{y\in\mathbb R}
e^{\lambda y}d\tilde \nu(y)
\leq \sup(\mathbb R\setminus [c_0,+\infty))=c_0\] 
holds. 
\hfill 
$\blacksquare$  

\sect{Proof of Theorem 2}
\( \, \, \, \, \, \, \, \) In this section, we prove Theorem 2. 
First, we give a basic fact for the linear equation
\begin{equation}v_t=\hat \mu*v\end{equation}
on the phase space $BC(\mathbb R)$: 
\begin{lemma} \ 
Let $\hat \mu$ be a Borel-measure on $\mathbb R$ with $\hat \mu (\mathbb R)<+\infty$. 
Let $\, \hat P:\, BC(\mathbb R)\rightarrow BC(\mathbb R)\, $ be the time $1$ map 
of the flow on $BC(\mathbb R)$ generated by the linear equation (6.1). 
Then, there exists a Borel-measure $\hat \nu$ on $\mathbb R$ 
with $\hat \nu(\mathbb R)<+\infty$ such that 
\[\hat P[v]=\hat \nu*v\] 
holds for all $v\in BC(\mathbb R)$. 
Further, if $v$ is a nonnegative, bounded and continuous function on $\mathbb R$, 
then the inequality 
\[v+\hat \mu*v\leq \hat \nu*v\]
holds. 
\end{lemma}
{\bf Proof.} \ 
Put a functional $\, \tilde P:\, BC(\mathbb R)\rightarrow \mathbb R\, $ as 
\[\tilde P[v]:=(\hat P[v])(0).\] 
Then, the functional $\tilde P$ is linear, bounded and positive. Hence, 
there exists a Borel-measure $\tilde \nu$ on $\mathbb R$ with 
$\tilde \nu(\mathbb R)<+\infty$ 
such that if a continuous function $v$ on $\mathbb R$ 
satisfies $\lim_{|x|\rightarrow\infty}v(x)=0$, then 
\begin{equation}\tilde P[v]=\int_{y\in \mathbb R}v(y)d\tilde \nu(y)\end{equation} 
holds. 

Let $v\in BC(\mathbb R)$. Then, there exists 
a sequence $\{v_n\}_{n=1}^\infty\subset BC(\mathbb R)$ 
with $\sup_{n\in\mathbb N,x\in\mathbb R}|v_n(x)|<+\infty$ 
and $\lim_{|x|\rightarrow\infty}v_n(x)=0$ for all $n\in\mathbb N$ 
such that $v_n\rightarrow v$ as $n\rightarrow\infty$ 
uniformly on every bounded interval. From Proposition 10, (6.2) 
and $\tilde \nu(\mathbb R)<+\infty$, we have 
\[\tilde P[v]=\lim_{n\rightarrow\infty}\tilde P[v_n]
=\lim_{n\rightarrow\infty}\int_{y\in \mathbb R}v_n(y)d\tilde \nu(y)
=\int_{y\in \mathbb R}v(y)d\tilde \nu(y).\] 
We take a Borel-measure  $\hat \nu$ on $\mathbb R$ 
with $\hat \nu(\mathbb R)<+\infty$ such that 
\[\hat \nu((-\infty,y))=\tilde \nu((-y,+\infty))\] 
holds for all $y\in \mathbb R$. 
Then, for any $v\in BC(\mathbb R)$, we have
\[(\hat P[v])(x)\equiv \tilde P[v(\cdot+x)]
\equiv\int_{y\in \mathbb R}v(y+x)d\tilde \nu(y)
\equiv(\hat\nu*v)(x).\] 

Let $v$ be a nonnegative, bounded and continuous function on $\mathbb R$. 
Then, in $t\in[0,+\infty)$, the function 
\[u(t,x):=v(x)+t(\hat\mu*v)(x)\]
is a sub-solution to (6.1), 
because of $v(x)\leq u(t,x)$. Hence, 
\[v+\hat \mu*v\leq \hat P[v]\]
holds. 
\hfill 
$\blacksquare$ 
\begin{lemma} \  
Let $\hat \mu$ be a Borel-measure on $\mathbb R$ with $\hat \mu (\mathbb R)<+\infty$.  
Suppose a constant $\gamma$ and a Lipschitz continuous function $g$ on $\mathbb R$ 
with $g(0)=0$ satisfy $\gamma<g^\prime(0)$. 
Let $\, \tilde P:\, BC(\mathbb R)\rightarrow BC(\mathbb R)\, $ be the time $1$ map 
of the flow on $BC(\mathbb R)$ generated by the linear equation 
\begin{equation}v_t=\hat\mu*v+\gamma v.\end{equation}
Let $\, \tilde P_0:\, BC(\mathbb R)\rightarrow BC(\mathbb R)\, $ be the time $1$ map 
of the flow on $BC(\mathbb R)$ generated by the equation 
\begin{equation}v_t=\hat\mu*v+g(v).\end{equation} 
Then, there exists a positive constant $\varepsilon$ such that the inequality  
\[\tilde P[v]\leq \tilde P_0[v]\] 
holds for all $v\in BC(\mathbb R)$ with $0\leq v\leq \varepsilon$.  
\end{lemma} 
{\bf Proof.} \ 
We take a positive constant $\varepsilon$ such that 
\begin{equation}h\in[0,(1+e^{\hat\mu(\mathbb R)+\gamma})\varepsilon] 
\ \Longrightarrow \ \gamma h\leq g(h)\end{equation} 
holds. Let a function $v\in BC(\mathbb R)$ satisfy $0\leq v\leq \varepsilon$. 
Then, we take the solution $\tilde v(t,x)$ to (6.3) with $\tilde v(0,x)=v(x)$. 
We see 
\[0\leq \tilde v(t,x)\leq e^{(\hat\mu(\mathbb R)+\gamma)t}\varepsilon
\leq (1+e^{\hat\mu(\mathbb R)+\gamma})\varepsilon\] 
for all $t\in [0,1]$. Hence, from (6.5), in $t\in [0,1]$, 
the function $\tilde v(t,x)$ is a sub-solution to (6.4). 
So, the inequality 
\[(\tilde P[v])(x)=\tilde v(1,x)\leq (\tilde P_0[v])(x)\]
holds. 
\hfill 
$\blacksquare$ 

\vspace*{0.8em} 

We use Proposition 15, Lemmas 16 and 17 to show the following: 
\begin{lemma} \  
Let $f^\prime(0)>0$. Suppose there exist $c\in\mathbb R$ 
and $\psi\in\mathcal M$ with $\psi(-\infty)=0$ and $\psi(+\infty)=1$ 
such that $\{\psi(x+ct)\}_{t\in\mathbb R}$ is a solution to (1.1). 
Then, there exists a positive constant $\lambda$ such that 
\[\int_{y\in\mathbb R}e^{-\lambda y}d\mu(y)<+\infty\] 
holds. 
\end{lemma} 
{\bf Proof.} \ 
Let $\hat \mu$ be the Borel-measure on $\mathbb R$ 
with $\hat \mu(\mathbb R)=1$ such that 
\[\hat \mu((-\infty,y))=\mu((-y,+\infty))\] 
holds for all $y\in \mathbb R$. 
Let $\, \hat P:\, BC(\mathbb R)\rightarrow BC(\mathbb R)\, $ be the time $1$ map 
of the flow on $BC(\mathbb R)$ generated by the linear equation (6.1). 
Then, by Lemma 16, there exists a Borel-measure $\hat \nu$ on $\mathbb R$ 
with $\hat \nu(\mathbb R)<+\infty$ such that for any $v\in BC(\mathbb R)$, 
\begin{equation}\hat P[v]=\hat \nu*v\end{equation} 
holds and for any nonnegative, bounded and continuous function $v$ on $\mathbb R$,  
\begin{equation}\hat \mu*v\leq \hat \nu*v\end{equation}
holds. Let $\, \tilde P:\, BC(\mathbb R)\rightarrow BC(\mathbb R)\, $ be the time $1$ map 
of the flow on $BC(\mathbb R)$ generated by the linear equation 
\[v_t=\hat \mu*v-v+\frac{f^\prime(0)}{2}v.\] 
Then, from (6.6) and (6.7), as $\tilde \nu$ is the Borel-measure on $\mathbb R$ 
defined by 
\[\tilde \nu:=e^{-1+\frac{f^\prime(0)}{2}}\hat \nu,\] 
we have 
\begin{equation}\tilde P[v]=\tilde \nu*v\end{equation} 
for all $v\in BC(\mathbb R)$ and   
\begin{equation}\hat \mu*v \leq e^{1-\frac{f^\prime(0)}{2}}(\tilde \nu*v)\end{equation}
for all nonnegative, bounded and continuous functions $v$ on $\mathbb R$.  
Because $\tilde \nu(\mathbb R)=(\tilde\nu*1)(0)
=(\tilde P[1])(0)=e^{\frac{f^\prime(0)}{2}}$ holds from (6.8), 
we also have 
\begin{equation}1<\tilde \nu(\mathbb R)<+\infty.\end{equation} 

Let $\, \tilde Q_0:\, \mathcal B\rightarrow \mathcal B\, $ be the time $1$ map 
of the semiflow on $\mathcal B$ generated by the equation (5.1). 
Then, from Lemma 17 and (6.8), 
there exists a positive constant $\varepsilon$ such that the inequality   
\[\tilde \nu*u=\tilde P[u]\leq \tilde Q_0[u]\] 
holds for all $u\in\mathcal B$ with $u\leq \varepsilon$. 
Further, $\hat\psi(x-ct):=\psi(-(x-ct))$ is a solution to (5.1). 
Therefore, by Proposition 15 and (6.10), we obtain the inequality 
\[\inf_{\lambda>0}\frac{1}{\lambda}\log\int_{y\in\mathbb R}
e^{\lambda y}d\tilde \nu(y)
\,
\leq 
\,
c.\] 
So, there exists a positive constant $\lambda$ such that 
\[\int_{y\in\mathbb R}e^{\lambda y}d\tilde \nu(y)
\leq e^{\lambda(c+1)}<+\infty\] 
holds. Hence, from (6.9), 
\[
\int_{y\in\mathbb R}e^{-\lambda y}d\mu(y)
=\int_{y\in\mathbb R}e^{\lambda y}d\hat\mu(y)
=\lim_{n\rightarrow\infty}\int_{y\in\mathbb R}
\min\{e^{\lambda y},n\}d\hat\mu(y)
\]
\[
=\lim_{n\rightarrow\infty}
(\hat\mu*\min\{e^{-\lambda x},n\})(0)
\leq e^{1-\frac{f^\prime(0)}{2}}
\lim_{n\rightarrow\infty}
(\tilde\nu*\min\{e^{-\lambda x},n\})(0)
\]
\[
=e^{1-\frac{f^\prime(0)}{2}}
\lim_{n\rightarrow\infty}\int_{y\in\mathbb R}
\min\{e^{\lambda y},n\}d\tilde \nu(y)
=e^{1-\frac{f^\prime(0)}{2}}\int_{y\in\mathbb R}
e^{\lambda y}d\tilde \nu(y)<+\infty
\]
holds. 
\hfill 
$\blacksquare$ 

\vspace*{0.8em} 

\noindent
{\bf Proof of Theorem 2.} 

It follows from Theorem 11 and Lemma 18. 
\hfill 
$\blacksquare$

\vspace*{1.6em} 

\noindent Acknowledgments. \ 
It was partially supported 
by Grant-in-Aid for Scientific Research (No.19740092) 
from Ministry of Education, Culture, Sports, Science 
and Technology, Japan. 

\newpage

\[ \begin{array}{c} \mbox{R\scriptsize EFERENCES}  \end{array} \]

[1] C. Atkinson and G. E. H. Reuter,  
Deterministic epidemic waves, 
{\it Math. Proc. Cambridge Philos. Soc.}, 80 (1976), 315-330. 

[2] P. W. Bates, P. C. Fife, X. Ren and X. Wang, Traveling waves in a convolution model 
for phase transitions, {\it Arch. Rational Mech. Anal.}, 138 (1997), 105-136. 

[3] J. Carr and A. Chmaj, Uniqueness of travelling waves for nonlocal monostable equations, 
{\it Proc. Amer. Math. Soc.}, 132 (2004), 2433-2439. 

[4] X. Chen, Existence, uniqueness, and asymptotic stability of traveling waves 
in nonlocal evolution equations, {\it Adv. Differential Equations}, 2 (1997), 125-160. 

[5] X. Chen and J.-S. Guo, Uniqueness and existence of traveling waves for discrete quasilinear monostable dynamics, {\it Math. Ann.}, 326 (2003), 123-146. 

[6] J. Coville, J. D\'{a}vila and S. Mart\'{i}nez, 
Nonlocal anisotropic dispersal with monostable nonlinearity, 
{\it J. Differential Equations}, 244 (2008), 3080-3118. 

[7] J. Coville and L. Dupaigne, On a non-local equation arising in population dynamics, 
{\it Proc. Roy. Soc. Edinburgh A}, 137 (2007), 727-755.


[8] R. A. Fisher, {\it The Genetical Theory of Natural Selection}, Clarendon Press, Oxford, 1930. 

[9] R. A. Fisher, The wave of advance of advantageous genes, {\it Ann. Eugenics}, 7 (1937), 335-369. 

[10] B. H. Gilding and R. Kersner, {\it Travelling Waves in Nonlinear 
Diffusion-Convection Reaction}, Birkh\"auser, Basel, 2004. 

[11] J.-S. Guo and F. Hamel, Front propagation for discrete periodic monostable equations, 
{\it Math. Ann.}, 335 (2006), 489-525. 

[12] J.-S. Guo and Y. Morita, Entire solutions of reaction-diffusion equations 
and an application to discrete diffusive equations, {\it Discrete Contin. Dynam. Systems}, 
12 (2005), 193-212. 

[13] F. Hamel and N. Nadirashvili, Entire solutions of the KPP equation, 
{\it Comm. Pure Appl. Math.}, 52 (1999), 1255-1276. 

[14] Y. Hosono, The minimal speed for a diffusive Lotka-Volterra model, {\it Bull. Math. Biol.}, 
60 (1998), 435-448. 

[15] Y. Kan-on, Fisher wave fronts for the Lotka-Volterra competition model with diffusion, 
{\it Nonlinear Anal.}, 28 (1997), 145-164. 

[16] A. N. Kolmogorov, I. G. Petrovsky and N. S. Piskunov, \'Etude de l'\'equation de la difusion 
avec croissance de la quantit\'e de mati\`ere et son application \`a un probl\`eme biologique, 
{\it Bull. Univ. Moskov. Ser. Internat. A}, 1 (1937), 1-25. 

[17] B. Li, H. F. Weinberger and M. A. Lewis, Spreading speeds as slowest wave speeds 
for cooperative systems, {\it Math. Biosci.}, 196 (2005), 82-98. 

[18] X. Liang, Y. Yi and X.-Q. Zhao, 
Spreading speeds and traveling waves for periodic evolution systems, 
{\it J. Differential Equations}, 231 (2006), 57-77.  

[19] X. Liang and X.-Q. Zhao, Asymptotic speeds of spread and traveling waves 
for monotone semiflows with applications, {\it Comm. Pure Appl. Math.}, 60 (2007), 1-40. 

[20] H. Okamoto and M. Shoji, {\it The Mathematical Theory of Permanent Progressive Water-Waves}, 
World Scientific Publishing Co., River Edge, 2001. 

[21] K. Schumacher, Travelling-front solutions for integro-differential equations. I, 
{\it J. Reine Angew. Math.}, 316 (1980), 54-70. 

[22] K. Schumacher, Travelling-front solutions for integrodifferential equations II, 
{\it Biological Growth and Spread}, pp. 296-309, Springer, Berlin-New York, 1980. 

[23] K. Uchiyama, The behavior of solutions of some nonlinear diffusion equations for large time, 
{\it J. Math. Kyoto Univ.}, 18 (1978), 453-508. 

[24] A. I. Volpert, V. A. Volpert and V. A. Volpert, {\it Traveling wave solutions 
of parabolic systems}, American Mathematical Society, Providence, 1994. 

[25] H. F. Weinberger, Long-time behavior of a class of biological models, 
{\it SIAM J. Math. Anal.}, 13 (1982), 353-396. 

[26] H. F. Weinberger, On spreading speeds and traveling waves for growth and migration models 
in a periodic habitat, {\it J. Math. Biol.}, 45 (2002), 511-548. 


[27] H. F. Weinberger, M. A. Lewis and B. Li, 
Anomalous spreading speeds of cooperative recursion systems,  
{\it J. Math. Biol.}, 55 (2007), 207-222. 

[28] H. Yagisita, Existence of traveling wave solutions 
for a nonlocal monostable equation: an abstract approach, 
{\it Discrete Contin. Dyn. Syst.}, submitted 
({\bf For editors and referees:} http://arxiv.org/abs/0807.3612). 

\end{document}